\documentclass[onecolumn, 12pt, conference]{ieeeconf}
\usepackage{geometry}                		
\usepackage[parfill]{parskip}    		
\usepackage{algorithmic,bm}
\usepackage{algorithm}
\usepackage{amssymb,url}
\usepackage{amsmath,amsfonts,mathrsfs}
\usepackage{graphicx}
\usepackage{caption,siunitx}
\usepackage{subfigure}
\usepackage[pagebackref=true,breaklinks=true,letterpaper=true,colorlinks,bookmarks=false,colorlinks]{hyperref}

\newcommand{\R}{\mathbb{R}}
\def\F{\vec{F}}

\newcommand{\Lt}{L^2(\Omega)}

\renewcommand{\u}{\vec{u}}

\newcommand{\curl}{\operatorname{curl}}

\graphicspath{{./figs/}}

\IEEEoverridecommandlockouts
\title{Where computer vision can aid physics: dynamic cloud motion forecasting from satellite images}
\author{Sergiy Zhuk\thanks{IBM Research, Dublin, Ireland, \texttt{\{sergiy.zhuk,tigran,albert\}@ie.ibm.com}} \and Tigran Tchrakian \and Albert Akhriev \and Siyuan Lu\thanks{IBM T.J. Watson Research Center, Yorktown, U.S., \texttt{\{lus,hendrikh\}@us.ibm.com}} \and Hendrik Hamann}
\begin{document}
\maketitle
\begin{abstract}
This paper describes a new algorithm for solar energy forecasting from a sequence of Cloud Optical Depth (COD) images. The algorithm is based on the following simple observation: the dynamics of clouds represented by COD images resembles the motion (transport) of a density in a fluid flow. This suggests that, to forecast the motion of COD images, it is sufficient to forecast the flow. The latter, in turn, can be accomplished by fitting a parametric model of the fluid flow to the COD images observed in the past. Namely, the learning phase of the algorithm is composed of the following steps: (i) given a sequence of COD images, the snapshots of the optical flow are estimated from two consecutive COD images; (ii) these snapshots are then assimilated into a Navier-Stokes Equation (NSE), i.e. an initial velocity field for NSE is selected so that the corresponding NSE' solution is as close as possible to the optical flow snapshots. The prediction phase consists of utilizing a linear transport equation, which describes the propagation of COD images in the fluid flow predicted by NSE, to estimate the future motion of the COD images. The algorithm has been tested on COD images provided by two geostationary operational environmental satellites from NOAA serving the west-hemisphere.
\end{abstract}
\section{Introduction and problem statement}
\label{sec:intro}
\textbf{Motivation.} Solar energy is the most abundant form of renewable energy resources and its contribution towards the total energy mix is rapidly increasing~\cite{US_DoE2012}. However, integration of high penetration level of solar energy in the electric grid poses significant challenge and cost~\cite{denholm_EP_2007}. The cloud movement, formation, dissipation and associated variable shading of solar panels may result in steep ramps of solar power being injected into the grid. The variability and uncertainty of solar power often forces the system operators to hold extra reserves of conventional power generation which adds cost. Ongoing research as well as previous experience of wind power integration shows that accurate solar forecasting plays a key role in the reliable and cost-effective integration of solar power~\cite{LuHamann_IEEE_Energy_2014}. For accurate short-term forecasting (a minute- to hour-ahead) of large geographical areas (continental scale), forecasting models using geostationary operational environmental satellite (GOES) imagery~\cite{Perez_SolarEnergy_2002} are more effective than numerical weather prediction models as the latter typically takes too long (several hours) to ramp up. The GOES satellite imagery only provides information on the current distribution of cloud, including cloud optical depth (COD) and top/bottom altitude which can be converted to the current solar irradiance at the earth’s surface using radiative transfer modeling. An accurate cloud advection model is thus required to forecast the future cloud distribution and solar irradiance.

\begin{figure}
        \centering
        \begin{subfigure}{}
                \includegraphics[width=3.5cm]{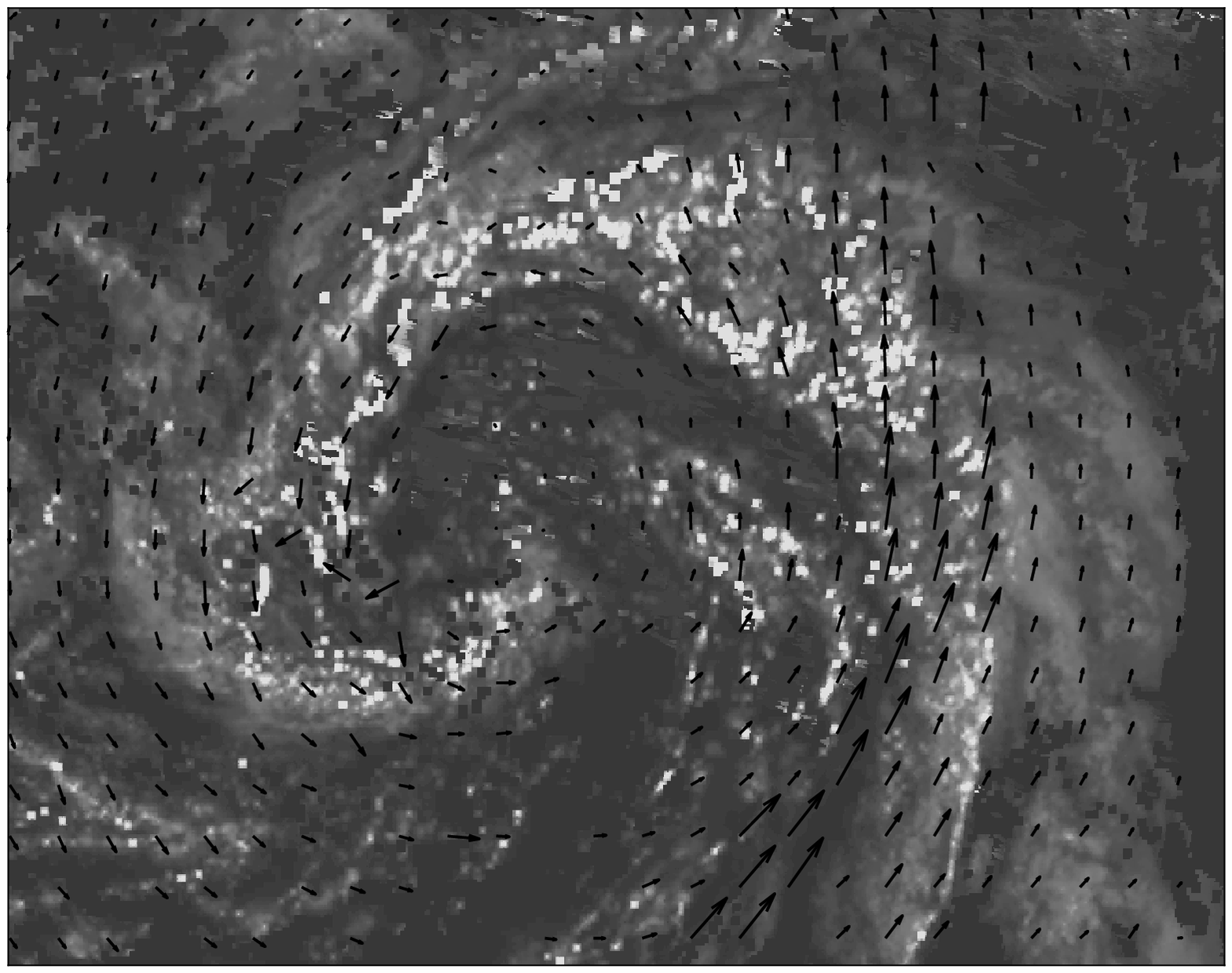}
        \end{subfigure}%
        \hspace{-3pt}
        \begin{subfigure}{}
                \includegraphics[width=3.5cm]{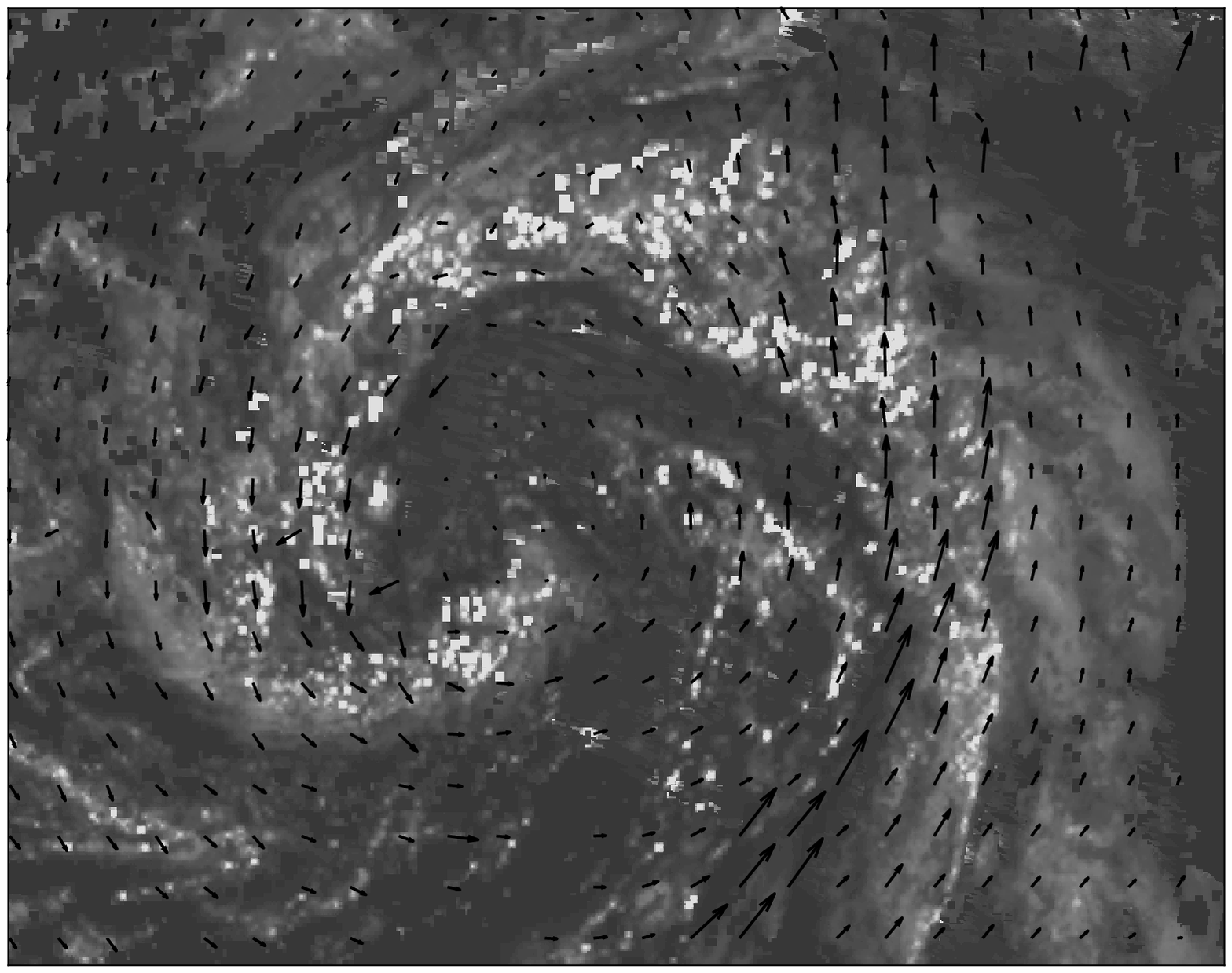}
        \end{subfigure}%
	\vspace{-3pt}
        \begin{subfigure}{}
                \includegraphics[width=3.5cm]{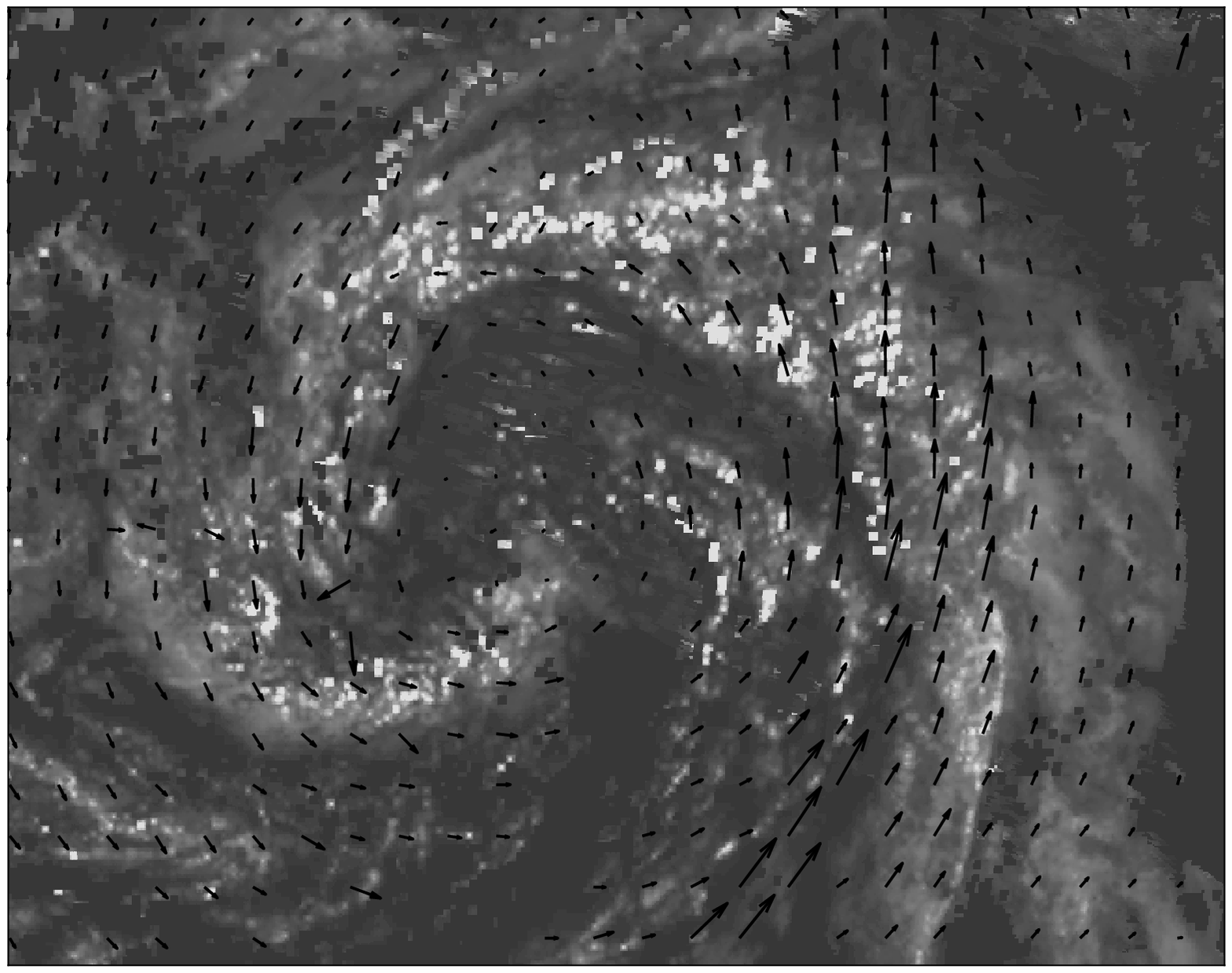}
        \end{subfigure}%
        \hspace{-6pt}
	 \begin{subfigure}{}
          \includegraphics[width=3.5cm]{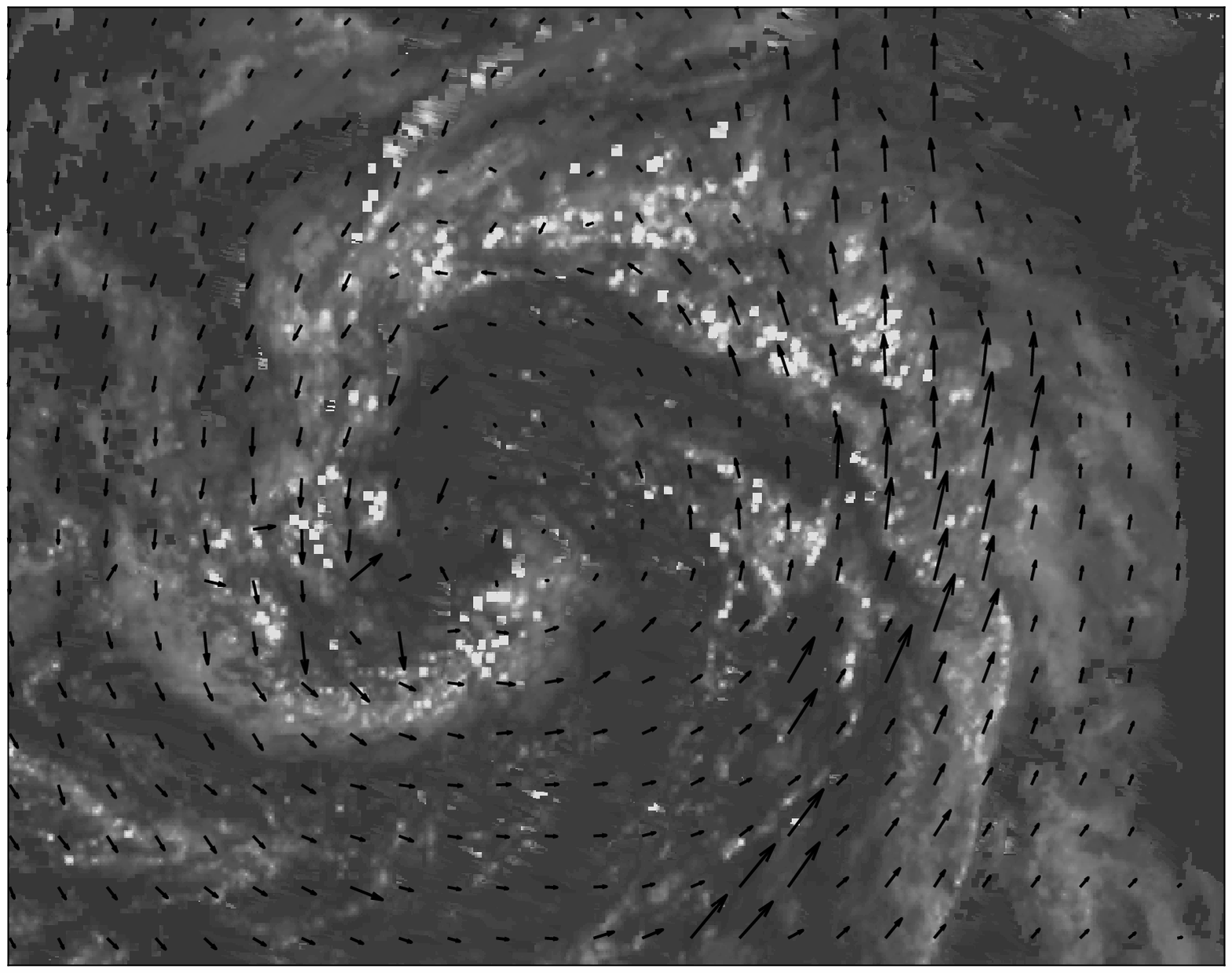}
        \end{subfigure}%
	\caption{Cloud optical depth and wind-field from the weather model on long-lat domain $[-140, -124]\times[39,51]$ at four consecutive times, 01.09.2013. A clear counter-clockwise rotation of the COD field is largely reflected in the wind-field. However, the wind-field suggests a strong upward flow towards the right of the domain which does not appear to be reflected in the COD field.}
\label{fig:wind}
\end{figure}

\textbf{Contribution.}  Conventionally cloud advections are performed using wind fields forecasted by numerical weather prediction (NWP) models~\cite{Miller2012} or wind fields determined using optical flow method \cite{Chow_SolarEnergy_2015}. In the former case, even when the NWP forecasted wind field is accurate, the inaccuracy in the estimate of the altitude of the clouds may still lead to significant error in the wind field and advected cloud distribution (see Fig.~\ref{fig:wind}). In the latter case, even though optical flow may accurately determine the current ``wind field'' at the altitude of the clouds, this  estimated ``wind field'' is assumed to persist for the forecasting period. In case when the actual wind field changes significantly in the hour ahead time scale, large error in the predicted cloud distribution may occur.

In this paper, we introduce a new algorithm for cloud motion prediction which assimilates purely data driven ``wind fields'' $\F_k=(u_k,v_k)^\top$ ($k$ -- discrete time index) derived by an optical flow estimation, into incompressible Navier-Stokes Equation (NSE) in two spatial dimensions. It builds upon a simple observation: the dynamics of clouds represented by COD images resembles the motion (transport) of a optical density $C(\bm{x},t)$, associated with a COD image, in a fluid flow $\u$. This suggests that, to predict the dynamics of the density $C(\bm{x},t)$ given a possibly sparse (in time) sequence of COD images, one can use the following procedure:
\begin{itemize}
\item [(i)] given a sequence of COD images, construct a map $\F_k=(u_k,v_k)^\top$, which transforms $C(\bm{x},t_{k-1})$ into $C(\bm{x},t_{k})$, by using some optical flow estimator;
\item [(ii)] assimilate $\F_k$ directly into the NSE, i.e. an initial velocity field for NSE is selected so that the corresponding NSE' solution, $\u$ passes through $\{\F_k\}$;
\item [(iii)] predict the dynamics of the velocity field $\u$ by integrating NSE forward in time until the desired horizon, and then utilize a linear transport equation (see equation~\eqref{eq:state} below) to predict the dynamics of $C(\bm{x},t)$.
\end{itemize}
To obtain maps $\F_k$ in step (i) we used an optical flow estimator from~\cite{sun13} which relies upon a combination of coarse-to-fine estimation by 3-level Gaussian pyramid, bicubic interpolation of image derivatives, and  robust penalty function augmented by median filtering.

\begin{figure}[h]
\centering
  \includegraphics[width=12cm]{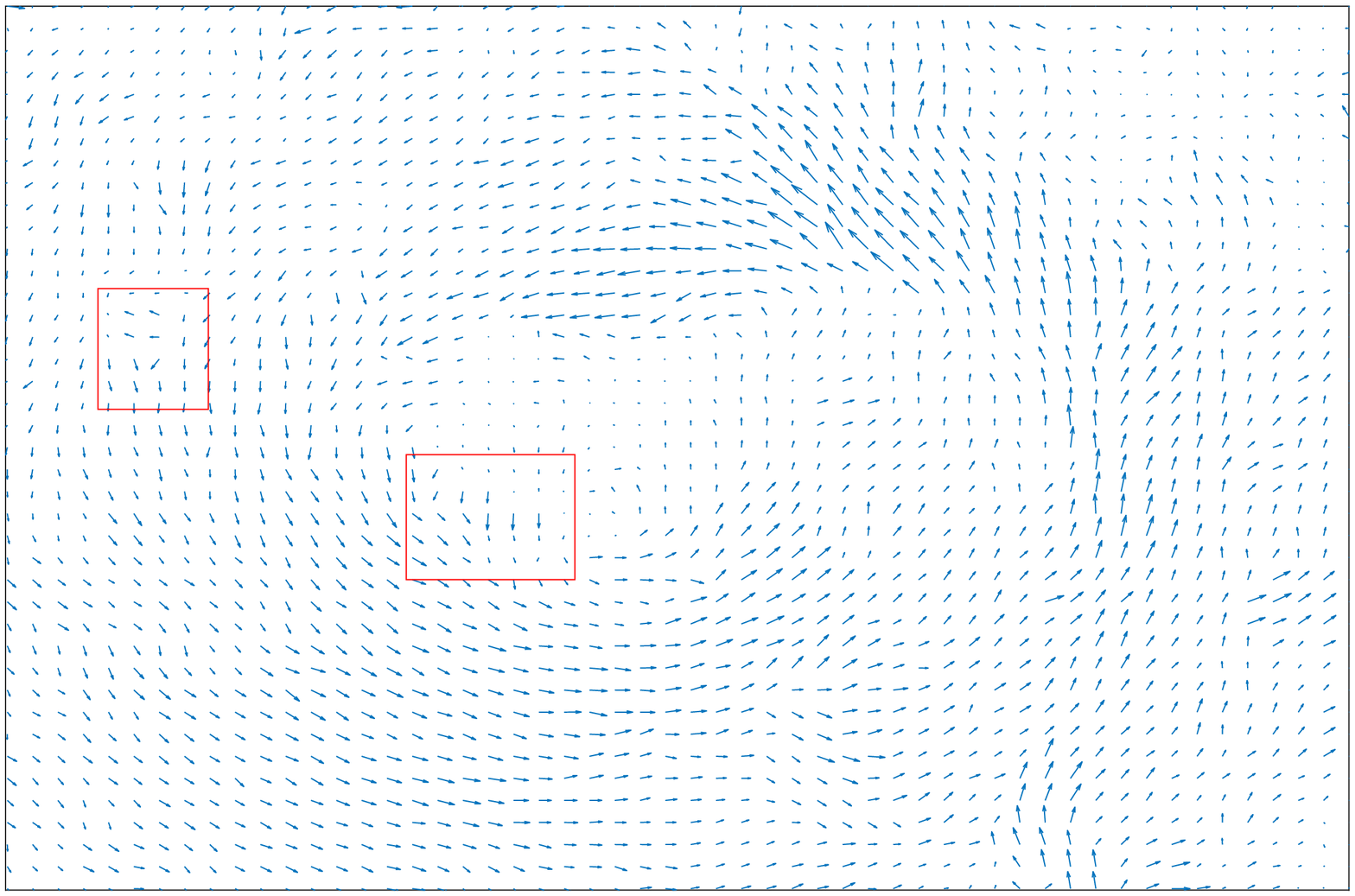}
  \caption{Optical flow estimate displaying a number of non-physical patterns (see red boxes).}
\label{fig:patches}
\end{figure}

Clearly, the quality of the estimated optical flow has strong impact onto the prediction accuracy. We stress that the COD images may be noisy due to the acquisition and transmission processes, and, in fact, this noise may mislead the estimator. In addition, it is well known that the optical flow estimation problem is generally ill-posed~\cite{HerlinBMZ12}, i.e. the COD images do not determine the unique map $\F_k$. Moreover, the COD images are usually sparse in time, and so they may not capture quite well fast underlying dynamic processes. Hence, the optical flow estimator has to rely upon a regulariser which may not account for the underlying physics of cloud motion. The latter may lead to some non-physical patterns in the resulting estimate of the optical flow (see Fig.~\ref{fig:patches}).

To overcome this issue, at step (ii) above we further regularize the data driven sequence of optical flows $\{\F_k\}$ obtained from (i): indeed, the solution of NSE captures some basic physics of the wind dynamics, and it introduces a spatio-temporal coherence between the subsequent optical flow estimates $\{\F_k\}$. Indirect assimilation of COD images into NSE has been previously considered in the literature (see~\cite{na11,HerlinBMZ12} and references therein). The approach taken here is different in that we assimilate the estimated optical flows directly into the NSE to simplify the computations and omit quite subtle numerical difficulties arising in the process of indirect assimilation. The precise differences are discussed section~\ref{s:problem_statement}. Finally, the prediction is done by advecting COD images forward in time by utilizing a discontinuous Galerkin (dG) approximation of the linear advection equation. As a result, the proposed algorithm enables accurate cloud  prediction at least 60 minutes ahead as it captures both the ``current wind field'' accurately determined by the data driven optical flow estimate, and some basic physics of the wind dynamics.

\textbf{Notation.} $\mathbb{C}^n/\mathbb{R}^n$ denotes the standard $n$-dimensional complex/real vector space, $x\cdot y:=\sum x_i\overline{y}_i$ denotes the canonical inner product in $\mathbb{C}^n/\mathbb{R}^n$, $\overline{g}$ denotes the complex conjugate of $g$, $A^\top$ denotes the transposed matrix, $A^\star$ denotes the conjugate transposed matrix, $\Omega\subset\R^2$ denotes a computational domain, $\Lt$ denotes the space of measurable squared-integrable functions over $\Omega$,\, $\|f\|^2_{\Lt}:=\int_\Omega f^2 dxdy$,\, $(f,g)_{\Lt}:=\int_\Omega f\overline{g} dxdy$,\, $\delta_{t-t_k}\psi(x,y,t):=\psi(x,y,t_k)$,\, $\Omega_T:=\Omega\times(t_0,T)$. $\partial_x$ and $\partial_y$ denote derivatives with respect to $x$ and $y$, $\curl\u=v_x-u_y$.

\subsubsection{Formal problem statement.}
\label{s:problem_statement}
Assume that COD images are observed over a domain of interest, $\Omega$ (in two spatial dimensions). We say that the cloud optical depth at position $\bm{x}:=(x,y)^\top\in\Omega$ and time $t$, $C(\bm{x},t)$ is transported by a velocity field $\u = (u(x,y,t),v(x,y,t))^\top$, which describes the flow of an incompressible fluid in $\Omega$, if
\begin{equation}\label{eq:state}
  \partial_t C + uC_x + vC_y=0\,.
\end{equation}
In computer vision literature, the latter equality is often reffered to as an optical flow constraint, and it has the following interpretation: each pixel $\bm{x}(t):=(x(t),y(t))^\top$ of a COD image moves with velocity $\u$, i.e. $\dot x(t) = u(x(t),y(t),t)$, $\dot y = v(x(t),y(t),t)$, and it carries the value of the optical density $C(\bm{x},t)$ along the trajectory of $\bm{x}$, i.e. $C(\bm{x}(0),0)=C(\bm{x}(t),t)$.

\begin{figure}[h]
  \centering
  \includegraphics[width=6cm]{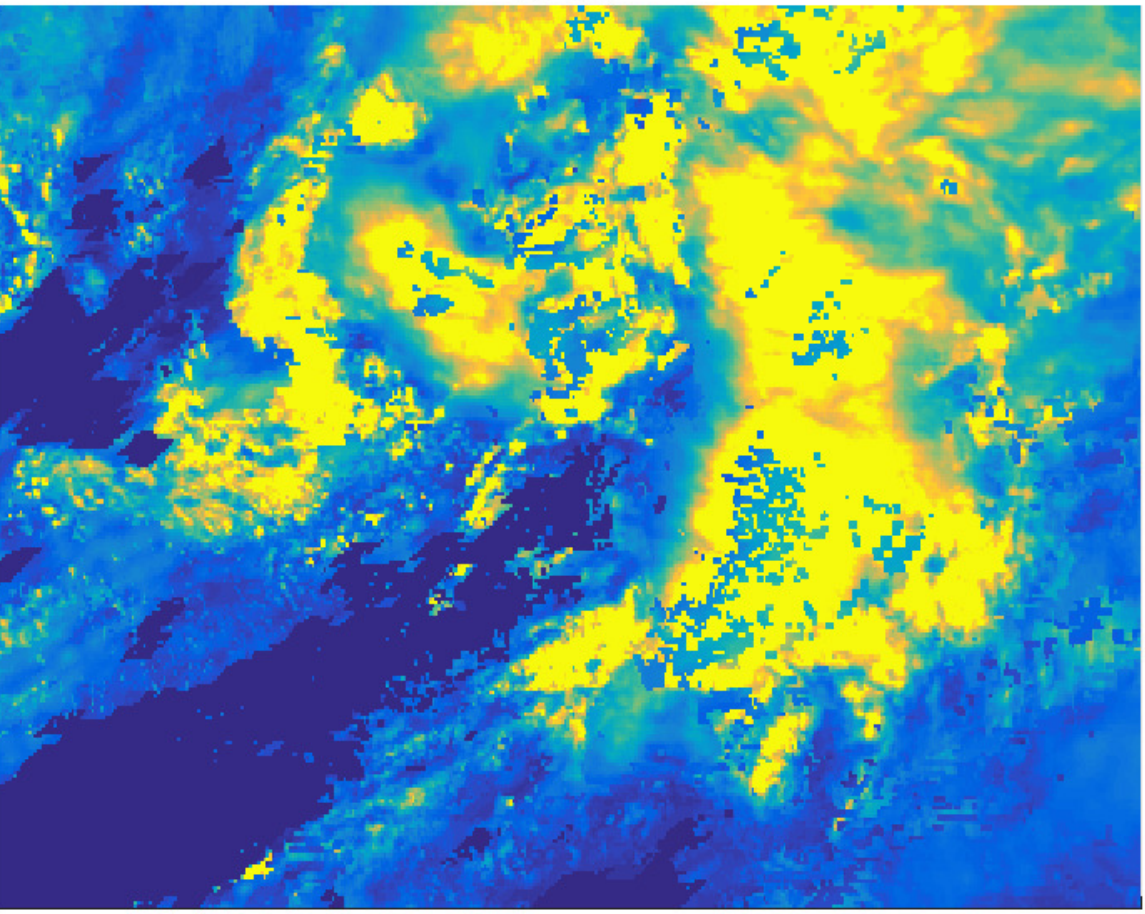}
  \includegraphics[width=6cm]{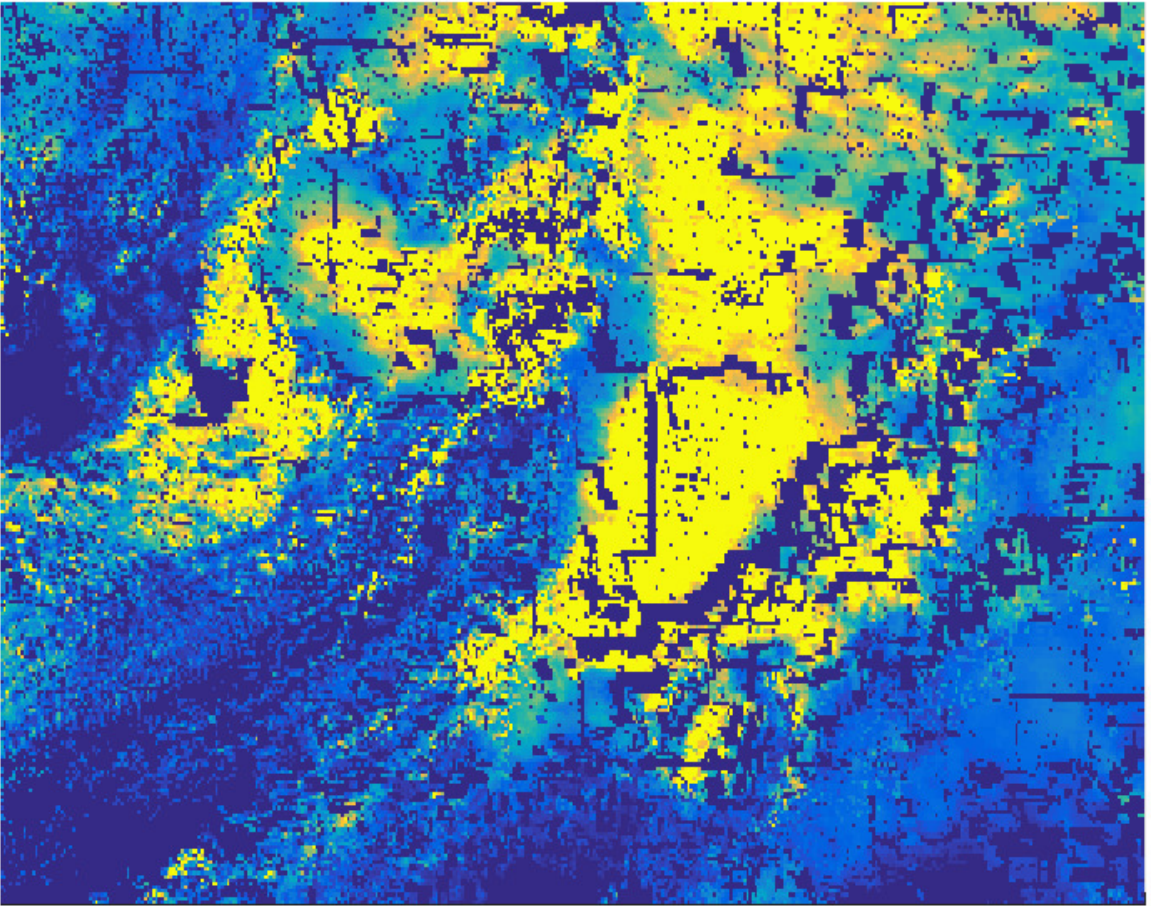}
  \caption{Advecting a COD image by propagating each COD pixel subject to the velocity field $\u$ at that pixel.\label{fig:groundtruth}}
\end{figure}

We stress that the straightforward advection (warping), i.e. propagating each pixel $\bm{x}(t)$ of a COD image subject to $\u$ at that pixel, may lead to developing artificial discontinuities as displayed on Figure~\ref{fig:groundtruth}.
The reason for this is that the warping is, in fact, equivalent to a very basic  numerical advection scheme, i.e. the position of the pixel $\bm{x}$ at time $t+1$ is obtained by using the forward Euler time integrator: $x(t+1) = x(t) + dt u(x(t),y(t),t)$, $y(t+1) = y(t) + dt v(x(t),y(t),t)$. In this work we use a more sophisticated discretization based on discontinuous Galerkin (dG) method which does not introduce artificial discontinuities.

\textbf{Dynamics of the fluid flow.} As noted above, the dynamics of COD images resembles the transport of a density/concentration of a mass in a fluid flow (in two spatial dimensions). This suggests to model $\u$ as the unique solution of the following two-dimensional Navier-Stokes equation in the vorticity-streamfunction weak formulation:
\begin{equation}
  \label{eq:NSE_periodic}
  \begin{split}
    &\dfrac{d}{dt}(\omega,\phi)_{\Lt} + b(\tilde u+\nabla^\perp\psi,\omega,\phi)+a(\omega,\phi)_{\Lt}
    =(f,\phi)_{\Lt}\,,\\ &a(\psi,\phi) = (\omega,\phi)_{\Lt}\\
    &\omega(0)=\curl(\u_0)\,, (x,t)\in\Omega_T
  \end{split}
\end{equation}
where $\phi$ is a test function, $\tilde u\in\R^2$ is the mean velocity field, $\omega$ is the vorticity, i.e. $\omega=\curl(\u)$, $\u_0$ is an initial velocity field, $f$ is a control parameter, e.g. source/sink, and \[
\begin{split}
b(\u,w,v)&:=(\u\cdot\nabla, w)_{\Lt}\,,a(\psi,\phi) = (\nabla\psi,\nabla\phi)_{\Lt}\\
& = (\psi_x,\phi_x)_{\Lt}+(\psi_y,\phi_y)_{\Lt}
\end{split}
\]
Note that the weak formulation above is obtained from the classical NSE by multiplying it with a test function $\phi$ in order to relax the smoothness assumptions on the solution $\omega$: in the above formulation, $\omega$ needs to have just first order derivatives.  We refer the reader to~\cite{MajdaBertozzi2002} for the further mathematical details on various fomulations of the NSE. Note that, by using integration by parts and selecting appropriate test/basis functions, it is not hard to show that the weak formulation~\eqref{eq:NSE_periodic} incapsulates various boundary conditions for the vorticity $\omega$, namely periodic, homogeneous Dirichlet and homogeneous Neumann. 

Note that $\u$ can be matched to a discrete (in time) sequence of COD images 
by recasting the classical optical flow estimation problem as an optimal control problem: select the control parameter $f$, and the initial velocity $\u_0$ such that the unique solution $C$ of~\eqref{eq:state}, corresponding to the unique solution $\u$ of~\eqref{eq:NSE_periodic}, is as close as possible (in some metric) to the given sequence of COD images (see~\cite{na11,HerlinBMZ12} and references therein). We stress that, for the case of COD images, the implementation of the latter control problem requires highly elaborate numerical schemes, and may be very expensive computationally. Indeed, transport model~\eqref{eq:state} is described by a hyperbolic equation which has Dirichlet boundary conditions on the inflow zone. The latter is a function of the coefficients $\u$ (see Section~\ref{sec:algorithm} for the details) which makes it quite complicated to design a numerically sound\footnote{i.e. the approximation which passes the so called gradient test} approximation of the gradient which is required to solve the aforementioned control problem numerically.

For this reason, we propose a hybrid data assimilation strategy: (i) construct a map $\F_k=(u_k,v_k)^\top$ which transforms COD image $C(\bm{x},t_{k-1})$ into the consecutive image $C(\bm{x},t_{k})$ by using an optical flow estimator, and (ii) to design an estimate of $\u$ which is close (locally) to each of the $\F_k$, and, at the same time, takes into account spatio-temporal patterns encoded in the COD images, and, importantly, is straightforward to implement. Now we are ready to formalize the problem statement:
\begin{itemize}
\item [A)] \emph{fitting problem:} given a sequence of optical flow maps $\{\F_k\}_{k=1}^{N_c-1}$ estimated from a sequence of COD images $\{C(\bm{x},t_k)\}_{k=1}^{N_c-1}$, find an initial vorticity field $\hat q$ and scalars $\hat u$, $\hat v$ such that the cost function
\begin{equation}
  \label{eq:OptContPr_1}
    J(\hat q,\hat u,\hat v):=\sum_{k=1}^{N_c} \|\left(
      \begin{smallmatrix}
        u_k - (\hat u +
    \delta_{t-t_k}\hat\psi_y)\\
    v_k -(\hat v -\delta_{t-t_k}\hat\psi_x)
      \end{smallmatrix}\right)\|^2_{\Lt}
  \end{equation}
is minimized along the solutions of~\eqref{eq:NSE_periodic}, provided $\tilde u=(\hat u,\hat v)^\top$.
\item [B)] \emph{advection/prediction problem:} given an estimate of $\u$ from $A)$ construct a short-term forecast of the future dynamics of $C$ by using~\eqref{eq:state}.
\end{itemize}
\section{Algorithm}
\label{sec:algorithm}

\subsubsection{Optical flow estimation.}
\label{sec:optflow}

As noted in the Introduction, we constructed maps $\F_k$ transforming $C(\bm{x},t_{k-1})$ into $C(\bm{x},t_{k})$ by using the optical flow estimator from~\cite{sun13} which aggregates state of the art optical flow estimaton techniques. This and a data preprocessing step are described below, in section~\ref{sec:experiments}.

\subsubsection{Fitting problem for NSE}
\label{sec:NSE}
Note that $J$ defined in~\eqref{eq:OptContPr_1} is a non-linear and non-convex functional. In what follows we provide a discrete representation of the gradient of $J$ and apply the L-BFGS version of the quasi-Newton method to solve the fitting problem (see point A) above).

To compute a discrete gradient of $J$ we proceed as follows. We set  $\Omega=(0,L_x)\times(0,L_y)$ and select the basis in $\Lt$ generated by complex-valued functions\footnote{In what follows, for a product of two basis functions, we always place the function of
$x$ to the left.} $\phi(x,y) = \phi_c(x) \phi_d(y)$, $\phi_n(x)=e^{2 \pi inz/L_x}$. This system is orthonormal in $\Lt$, i.e. $(\phi_i\phi_j,\phi_c\phi_d)_{\Lt}=\delta_{i,c}\delta_{j,d}$, provided the inner product is rescaled as follows: $(f,g)_{\Lt} := \frac{1}{L_xL_y}(f,g)_{\Lt}$. Now, denote by $\hat \omega$ and $\hat \psi$ the solution of~\eqref{eq:NSE_periodic} which corresponds to the initial vorticity $\hat\omega(0)=\hat q$ and $f=0$. We approximate the gradient of $J$ by using Fourier-Galerkin (FG) method. To this end, we define the projection operator $P_{N_x,N_y}$: \[
P_{N_x,N_y} v:=\sum_{\substack{
   |n|\leq N_x/2 \\
   |m|\leq N_y/2
  }} v_{nm}(t) \phi_n(x) \phi_m(y)\,.
\]
It projects any $v\in\Lt$ onto a $(N_x+1)(N_y+1)$-subspace of $L^2(\Omega)=L^2(0,L_x)\times L^2(0,L_y)$, which is generated by functions $\Phi:=\{\phi_c\phi_d\}_{|c|\le\frac{N_x}2,|d|\le\frac{N_y}2}$. Now, Fourier-Galerkin method suggests (i) to substitute $P_{N_x,N_y} \hat q$, $P_{N_x,N_y} \omega $ and $P_{N_x,N_y}\psi$ into~\eqref{eq:NSE_periodic}, and (ii) to plug the test function $\phi = \phi_c\phi_d$, for any $|c|\le\frac{N_x}2$ and $|d|\le\frac{N_y}2$, into the resulting system to generate a system of Ordinary Differential Equations (ODE) which defines the dynamics of the vector of the projection coefficients $\hat{\bm\omega}:=\{\hat\omega_{pq}\}_{|p|\le\frac{N_x}2,|q|\le\frac{N_y}2}$. This ODE takes the following form:
\begin{equation}
  \label{eq:FG_model}
  \dfrac{d\hat{\bm \omega}}{dt} = -B(\hat{\bm \omega})\hat{\bm\omega} - B(\tilde u) \hat{\bm\omega} - A\hat{\bm \omega}\,, \hat{\bm\omega}(0)=\hat{\bm q}\,,
\end{equation}
where $
B(\hat{\bm{\omega}}) = \{\sum_{p,q}\frac{\hat{\omega}_{pq}(pm-qn)L_xL_y}{p^2 L_y^2+q^2 L_x^2}\delta_{p+n,c} \delta_{q+m,d}\}$ with $|c|,|n|\le \frac{N_x}2$ and $|m|,|d|\le \frac{N_y}2$ is a skew-symmetric matrix approximating the non-linear advection term in~\eqref{eq:NSE_periodic}, and $
A:=\operatorname{diag}(\lambda_{-\frac{N_x}2,-\frac{N_y}2}\dots \lambda_{\frac{N_x}2,\frac{N_y}2})$ is a symmetric non-negative definite matrix composed of eigenvalues of the Laplacian, $\lambda_{cd}:=\frac{4\pi^2 c^2}{L_x^2} + \frac{4\pi^2 d^2}{L_y^2}$, $c^2+d^2>0$, $\lambda_{0,0}=0$. We refer the reader to~\cite{MajdaBertozzi2002} for the further details on FG method. 

Now, given $\hat q$ we can approximate $J$ as follows: $J(\hat q,\hat u,\hat v)\approx J(P_{N_x,N_y}\hat q,\hat u,\hat v)$. As a result, $J$ depends on a finite number of parameters: the vector of projection coefficients $\hat{\bm q}$ representing $P_{N_x,N_y}\hat q$, $\hat u$ and $\hat v$. Effectively, to each trajectory of~\eqref{eq:FG_model} which corresponds to the initial condition $\hat{\bm q}$ and coefficients $\hat u$ and $\hat v$, $J$ assignes a number, the Euclidian distance between the velocity field $u(x,y,t_k)=\hat\psi_y(x,y,t_k)$, $v(x,y,t_k) = \hat\psi_x(x,y,t_k)$, computed from the trajectory, and the corresponding optical flow map $\F_k=(u_k,v_k)^\top$. In other words, minimizing $J$ is equivalent to projecting $\F_k$ onto the manifold generated by trajectories of~\eqref{eq:FG_model}.\\
To fit $\hat{\bm q}$, $\hat u$, $\hat v$ to the ``data'' $\F_k=(u_k,v_k)^\top$ we need to compute the gradient of $J$ w.r.t. to $\hat{\bm q}$, $\hat u$, $\hat v$. The latter can be done by using a so-called adjoint method: one computes Gateaux derivative of $J$ w.r.t. $\hat{\bm q}$, obtaines a linear ODE for the gradient, and, finally, uses its adjoint equation to find the gradient explicitly. Here, we denote the solution of the adjoint equation by $\bm{\lambda}_k$ and form the gradient as follows: $\nabla J(\hat{\bm q},\hat u,\hat v)=2\sum_{k=1}^{N_c}\Psi_{1,\imath}^\star\bm{\lambda}_k(0),J_{\hat u},J_{\hat v})^\top$. $\bm{\lambda}_k$, $\Psi_{1,\imath}$, $J_{\hat u}$ and $J_{\hat v}$ are defined in the appendix~\ref{sec:gradJ}. This gradient is then used to approximate the solution of the fitting problem~\eqref{eq:OptContPr_1} by means of the L-BFGS method~\cite{L_BFGS} provided $\F_0=(u_0,v_0)^\top$ is used as the starting point.

\subsubsection{Advection/prediction problem.}
\label{sec:dG}
In what follows we apply the dG method for solving a linear advection equation with non-stationary coefficients:
\begin{equation}
\label{eq:advection2}
\partial_t C(\bm{x},t)+ \u \cdot \nabla C(\bm{x},t)=0, \quad \bm{x} \in \Omega \in \mathbb{R}^2
\end{equation}
with Dirichlet boundary conditions $C(\bm{x},t)=g$ on the inflow zone of the boundary of $\Omega$, $\partial\Omega$. The inflow zone is composed of all $\bm{x}\in \partial\Omega$ such that $\u(\bm{x},t)$ points inside $\Omega$. Since $\u$ varies in time, the inflow zone could ``travel'' along the boundary of $\Omega$ reflecting changes of $\u$. This makes it quite complicated to fit the initial vorticity $\hat q$ to the COD images. Indeed, $C$ depends on $\u$ which depends on $\hat q$. Hence, to find $\hat q$ by minimizing some distance between $C$ and the COD images, one would need to differentiate the weak formulation~\eqref{eq:weakStatement} of~\eqref{eq:advection2} w.r.t. the inflow zone. The latter does not seem to be practical from the computational stand-point. One could have overcome this by resorting to the vanishing viscosity method but this approach would require one to deal with numerical boundary layers. To overcome these difficulties, we are assimilating $\F_k$ into~\eqref{eq:NSE_periodic} directly, and then use \eqref{eq:advection2} just for the advection/prediction.

Following~\cite{Hesthaven_dGNodal_2008} we approximate $C$ by using dG approach as follows. $\Omega$ is divided into $K$ of non-overlapping rectangular elements, $D^k$, i.e., $\Omega \simeq \Omega_h = \bigcup\limits_{k=1}^{K} D^k$. On each element $D^k$, $C$ is approximated by $C_h^k$, which is expressed as the series,
\begin{equation}
\label{eq:seriesDG}
C_h^k(\bm{x},t)= \sum_{i=1}^{N+1} C^k(\bm{x}_i^k,t) \ell_i^k(\bm{x}), \quad \bm{x} \in D^k,
\end{equation}
where $\ell_i^k(\bm{x})$ are Lagrange interpolating polynomials in two dimensions defined by Legendre-Gauss-Lobatto (LGL) points $\bm{x}_i$ (see~\cite{Hesthaven_dGNodal_2008}). Substituting $C_h^k$ into~\eqref{eq:advection2}, we form the residual $R_h^k$ on a single element:
\begin{equation}
\label{eq:residual}
R_h^k = \partial_t C_h^k(\bm{x},t)+ \u \cdot \nabla C_h^k(\bm{x},t), \quad \bm{x} \in D^k.
\end{equation}
The Galerkin methods involves taking $\ell_i^k$ as test functions (i.e. same as the expansion/trial functions) and forcing the residual to be orthogonal to each of these test functions. Doing this, and then using integration by parts to move the spatial derivatives off the state and onto the test functions gives the following weak statement on element $D_k$:
\begin{equation}
\label{eq:weakStatement}
\begin{split}
\int_{D^k}  \left( C_h^k(\bm{x},t)_t \ell_n^k(\bm{x}) - \bm{f}_h^k(\bm{x},t) \cdot \nabla \ell_n^k(\bm{x}) \right)  d\bm{x} = \\
-\int_{\partial D^k} \bm{\hat{n}}\cdot \bm{f}^* \ell_n^k(\bm{x}) d\bm{x} , \quad n=1 \ldots N+1,
\end{split}
\end{equation}
where $\bm{\hat{n}}$ is the outward facing unit normal, $\bm{f}_h^k=(uC_h^k,vC_h^k)^\top$ and $\bm{f}^*$ is the numerical flux function which we take to be the local Lax-Friedrichs flux:
\begin{equation}
\label{eq:LF}
\bm{f}^*(C_i,C_e,\bm{u}_i,\bm{u}_e) = \frac{C_i \bm{u}_i + C_e \bm{u}_e }{2}+\frac{c_s}{2}\bm{\hat{n}}(\bm{u}_i-\bm{u}_e),
\end{equation}
where subscripts $i$ and $e$ refer respectively to the interior and exterior values at a point on the boundary, and $c_s$ is the maximum absolute value of the signal speed normal to the boundary at that point, i.e.,
\begin{equation}
\label{eq:LFconst}
c_s=\max_{j\in\{i,e\}} | \hat n_x u_j + \hat n_y v_j |.
\end{equation}
The surface integral in~\eqref{eq:weakStatement} allows the elements to `communicate' with one another. Since we are using rectangular elements, this surface integral is the sum of four line integrals, each one over one face of the element. The exterior solution values in \eqref{eq:LF} and \eqref{eq:LFconst} refer to the values of $C$ on a neighbouring element sharing the node in question. In the case that this node lies on the physical boundary, $\partial \Omega$, the exterior value is determined by a ``physical'' Dirichlet boundary condition $g$. If flow direction at that node is `into' the domain, $\Omega$, then the exterior solution values is set to a prescribed value subject to the boundary conditions (inflow zone). In this case, $\bm{\hat{n}}\cdot \u<0$ at a point on the boundary. If, on the other hand, the flow direction at that boundary point is `out of' the domain (i.e. $\bm{\hat{n}}\cdot \bm{u}>0$), then a free exit boundary condition is imposed at that point by setting the exterior values equal to the interior values.

The weak DG formulation on a single element $D^k$ can be written as:
\begin{equation}
\label{eq:weakDG}
M^k \dfrac{d}{dt} \bm{I}_h^k - S_x^\top \bm{F}_x - S_y^\top \bm{F}_y =-\sum\limits_{i=1}^4 (-1)^i M_e^{k,i} \bm{F}^*_i
\end{equation}
where $M_k$, $S_x$ and $S_y$ are the mass and stiffness matrices with the latter corresponding to advection in the $x$- and $y$-directions. The vectors $\bm{I}_h^k$ is a grid-function representing $C_h^k$ on the $(N+1)^2$ quadrature points on element $D^k$, and the vectors $\bm{F}_x$ and $\bm{F}_y$ are grid functions representing the first and second components of $\bm{f}_h^k$ respectively. The matrices, $M_e^{k,i}$ are \emph{edge}-mass matrices for element $D^k$ on face $i$, where the faces, $i= 1 \ldots 4$, are ordered: \emph{left}, \emph{right}, \emph{lower}, \emph{upper}. These matrices act on the vector $\bm{F}^*_i$ which represents the numerical flux over each node on face $i$. An easy way to define the mass and stiffness matrices for 2D rectangular elements is to do so in terms of LGL quadrature points and weights in the 1D interval $\mbox{I}=[-1,1]$ and also in terms of Lagrange polynomials defined on this interval. We refer the reader to~\cite{SZTT_SISC17} for the further implementation details.

\section{Experiments}
\label{sec:experiments}

\begin{figure}[h]
  \centering
  \includegraphics[width=9cm]{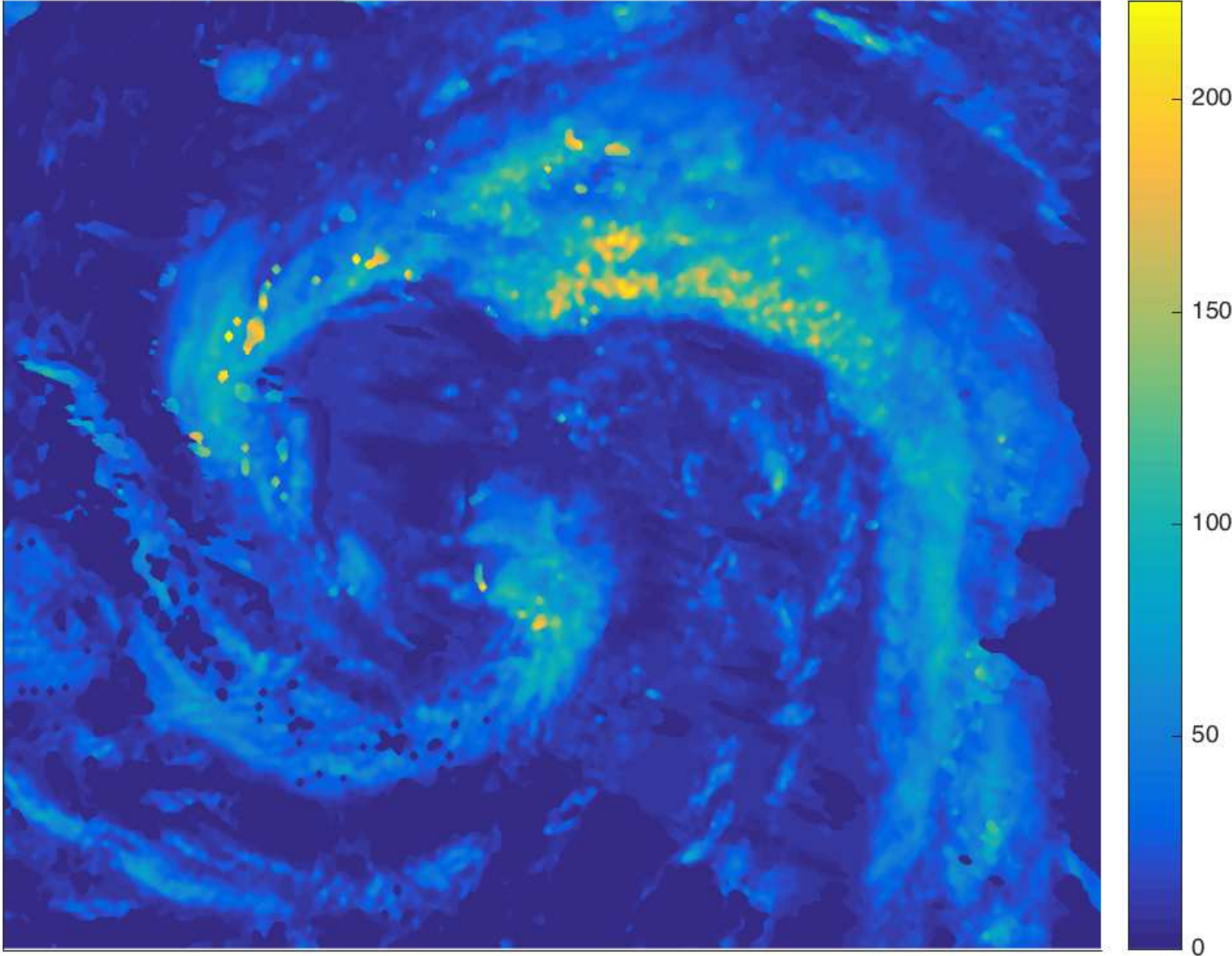}
  \caption{Example of a COD image over $\Omega$, 01.09.2013, 21:15. \label{fig:cod0}}
\end{figure}

\paragraph{COD images} The aforementioned algorithm has been tested on data from two GOES serving the west-hemisphere, GOES-15, the WEST satellite located over the Equator at \ang{135} west longitude and GOES-13, the EAST satellite located over the Equator at \ang{75} west longitude. The COD products of 15 minute interval and ~4 km spatial resolution are derived from GOES data using CLAVR-x~\cite{Thomas_Met_2004} made available by the NOAA Advanced Satellite Products Branch and Cooperative Institute for Meteorological Satellites Studies at the University of Wisconsin – Madison’s Space Science and Engineering Center (SSEC).
All data are on a longitude-latitude grid with uniform \ang{0.02} spacing. The area $\Omega$ covered by the data extends from \ang{-140} East (longitude) and \ang{39} North (latitude) to \ang{-124} degrees East (longitude) and \ang{51} North (latitude). For the experiment we used COD images recorded at 1 September 2013 over $\Omega$. (see Figure~\ref{fig:cod0}).

\begin{figure}
        \centering
        \begin{subfigure}{}
                \includegraphics[width=5.6cm]{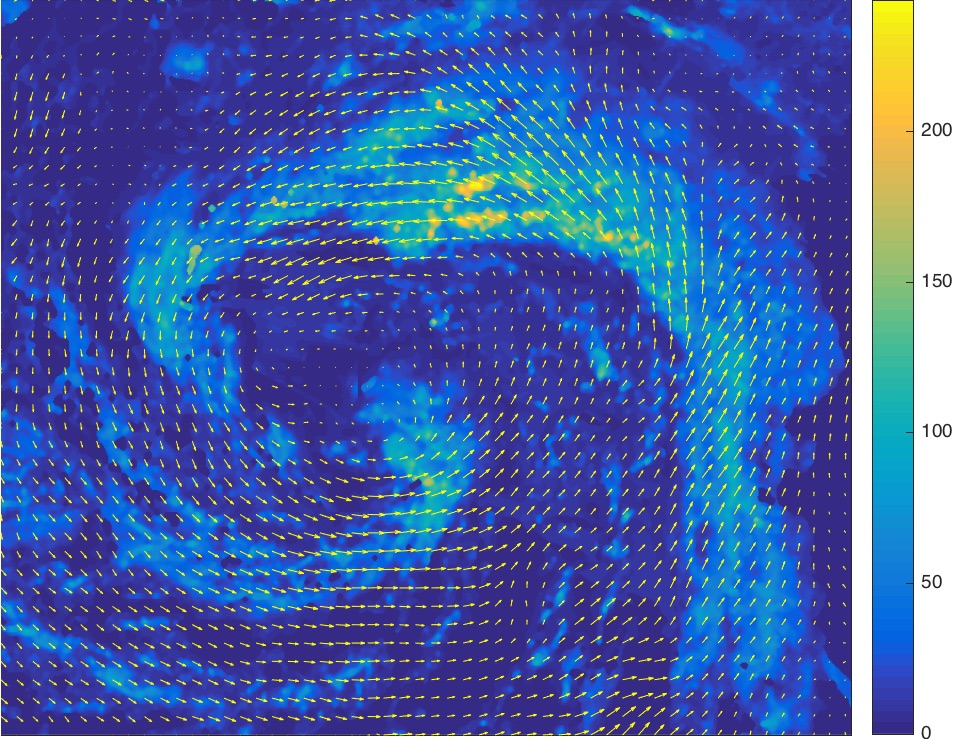}
        \end{subfigure}%
        \hspace{-7pt}
        \begin{subfigure}{}
                \includegraphics[width=5.6cm]{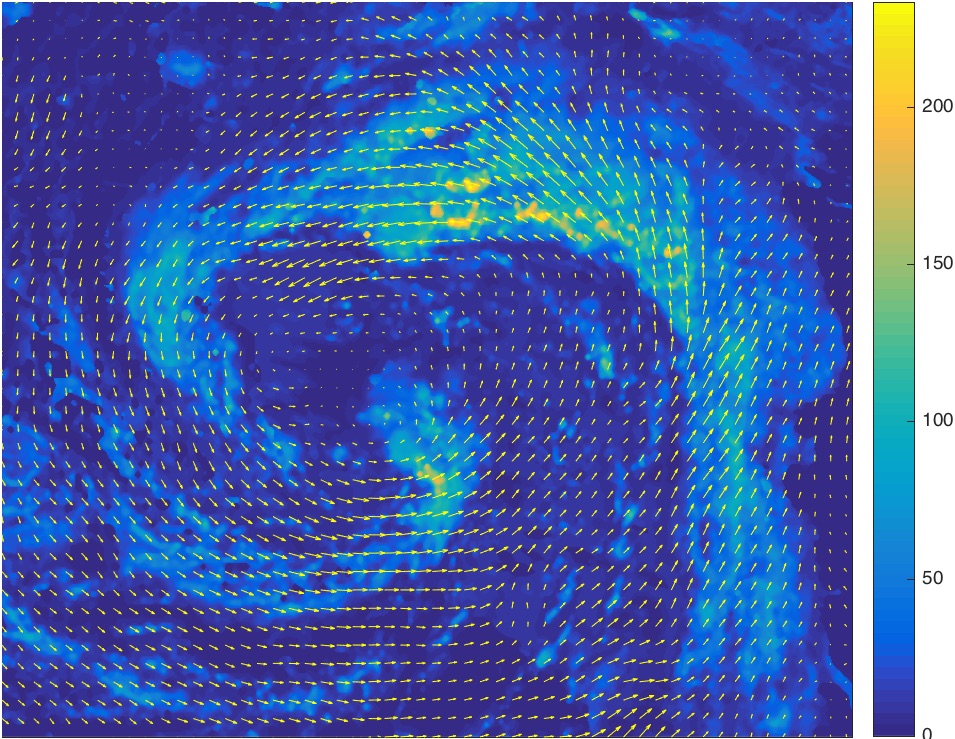}
        \end{subfigure}%
	\vspace{-2pt}
        \begin{subfigure}{}
                \includegraphics[width=5.6cm]{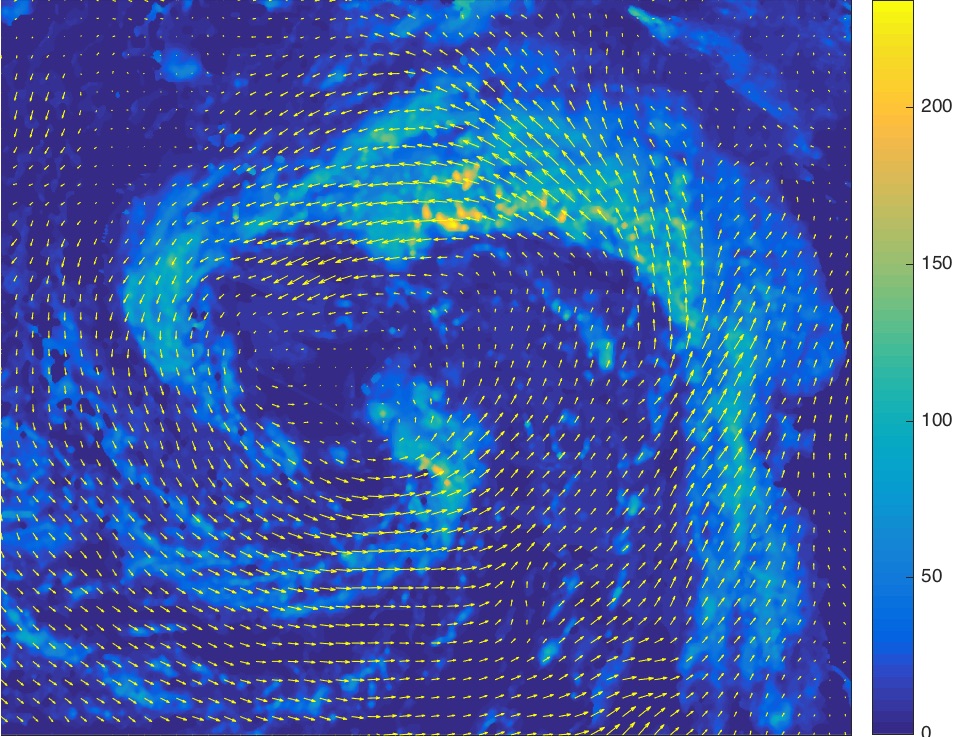}
        \end{subfigure}%
        \hspace{-10pt}
	 \begin{subfigure}{}
          \includegraphics[width=5.6cm]{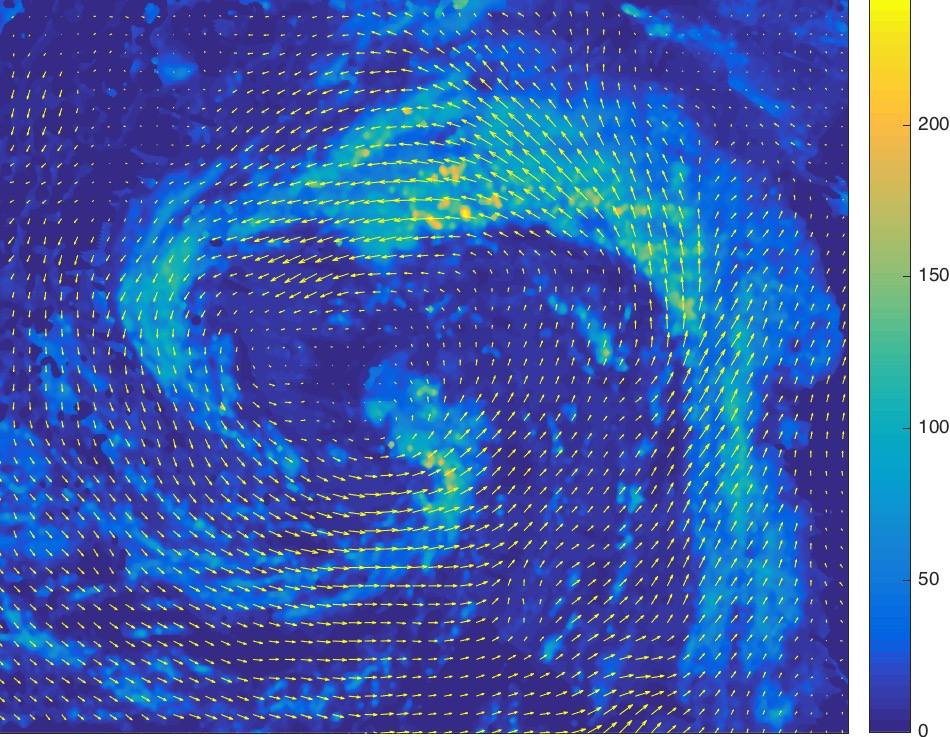}
        \end{subfigure}%
	\caption{COD images/optical flows $\F_k$ (yellow arrows) at four consecutive times, 21:30 - 22:15, 01.09.2013.}
\label{fig:cods}
\end{figure}

\paragraph{Optical flow estimates} We estimated the optical flow maps $\F_{1}\dots \F_{4}$ from the COD images displayed in Figure~\ref{fig:cods} by using the publicly available Matlab implementation of the optical flow estimator from~\cite{sun13}. It is important to note that the COD images required some preprocessing prior to invoking the estimator: without preprocessing the resulting maps $\F_k$ were of little or no use. The preprocessing step involved a number of transformations. First, the intensities of bright regions of a COD image are usually very high compared to the average intensity outside these regions. Careful data inspection revealed that the distortion preserves the local ordering and the local pattern. The following \textit{log}-transform was used to alleviate the distortion of intensity at each image pixel $(x,y)$: $\tilde{C}(x,y) =
\sigma \cdot \sqrt{\log{(1 + (C(x,y)/\sigma)^2)}/\log{(2)}}$, where $\sigma$ is a threshold value. If $C(x,y) = \sigma$, then $\tilde{C}(x,y) = C(x,y)$. For $C(x,y)\,{<}\,\sigma$ the $\tilde{C}(x,y)$ is nearly a linear function of $C(x,y)$, whereas for $C(x,y)\,{>}\,\sigma$ the logarithm dominates in the transformation. The threshold value $\sigma$ was selected as the $90$th percentile, which indicates the magnitude such that $90$\% percent of the pixel intensities have values less than this number. In other words, $10$\% of intensities were considered as potential outliers and those were suppressed by the logarithm. The optical flow was computed on transformed images with noticeable improvement comparing to the untransformed case.

Finally, the mean brightness of consecutive images $C(\bm{x},t_{k})$ and $C(\bm{x},t_{k+1})$ tends to variate significantly possibly due to some problems occurring in the acquisition/transmission processes. This variation was compensated by equating the median intensity values of $C(\bm{x},t_{k})$ and $C(\bm{x},t_{k+1})$: the medians are computed over the non-zero (``day'') pixels and the first image $C(\bm{x},t_{k})$ is scaled to make both medians equal.

\begin{figure}[h]
  \centering
  \includegraphics[width=9cm]{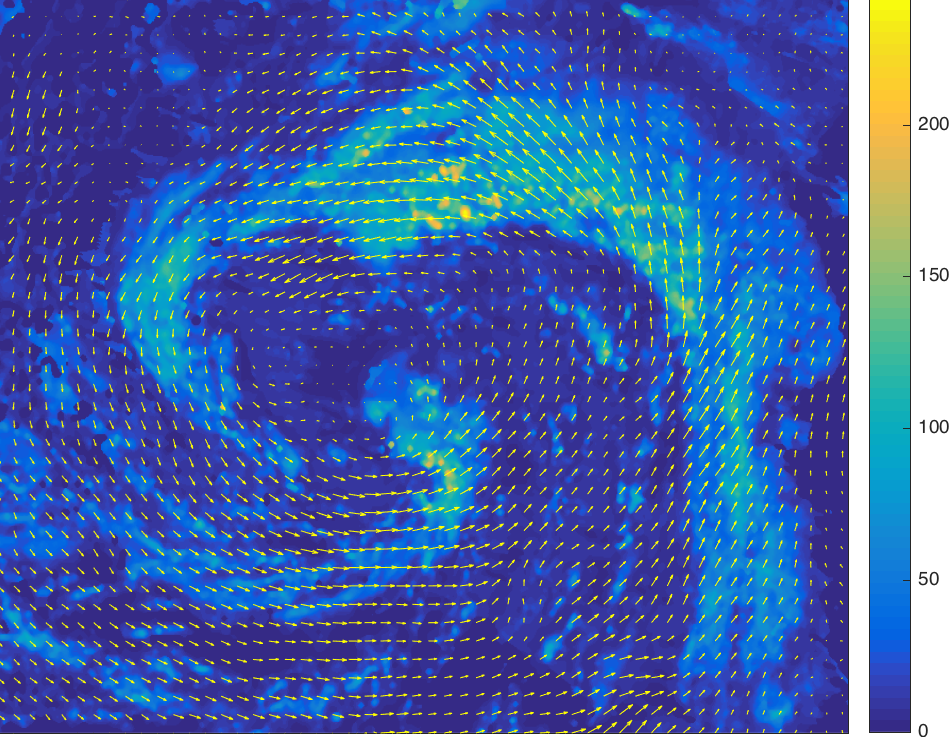}
  \caption{COD image/optical flow map $\F_4$ at (yellow arrows) 22:15\label{fig:cod1F4}}
\end{figure}

After the data preprocessing step, $\F_1$ was estimated from two consecutive (15 minutes apart) images (Fig.~\ref{fig:cods}, top left image, 21:30) and so on. The resulting maps are presented in Figure~\ref{fig:cods},\ref{fig:cod1F4}.

\begin{figure}[h]
  \centering
  \includegraphics[width=12cm]{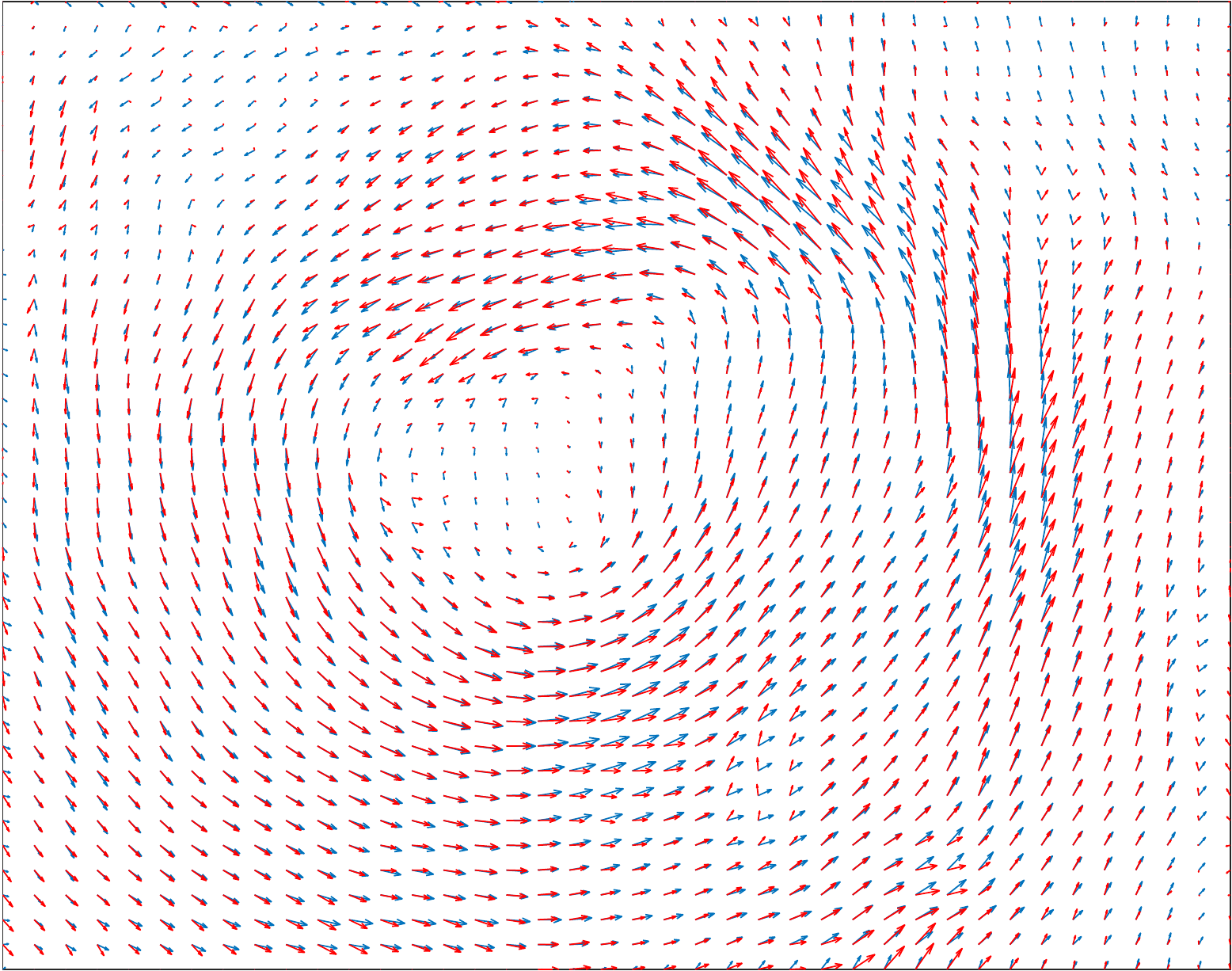}
  \caption{Optical flow $\F_4$ (red) at 22:15 versus NSE velocity field (blue) at 22:15\label{fig:F4NSE}}
\end{figure}

To fit NSE solution to the COD images we used the algorithm of Section~\ref{sec:NSE}. We took $N_x=N_y=28$, $L_x=1779$ and $L_y=941.8$ and timestep was taken to be the same as in the dG advection (see below). The NSE velocity field is compared versus $\F_4$ in Figure~\ref{fig:F4NSE}. As expected, the NSE velocity looks more regular.

\begin{figure}[!ht]
  	\centering
    	\includegraphics[width=2.5cm]{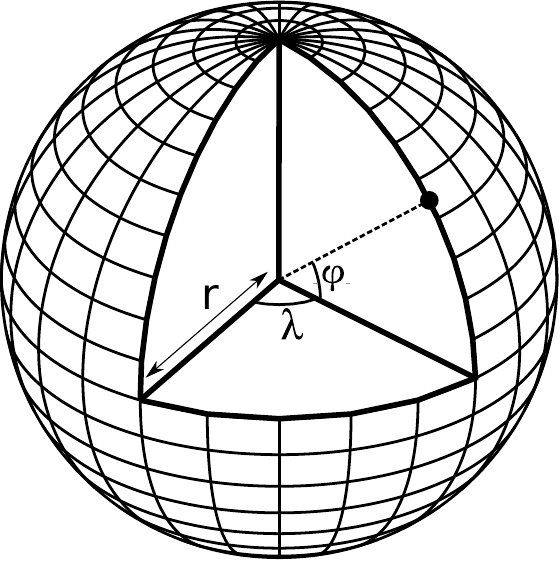}
	\caption{Longitude-latitude ($\lambda$-$\varphi$) coordinate system.}
	\label{fig:latLongGlobe}
\end{figure}

\begin{figure}[!ht]
  	\centering
    	\includegraphics[width=6cm]{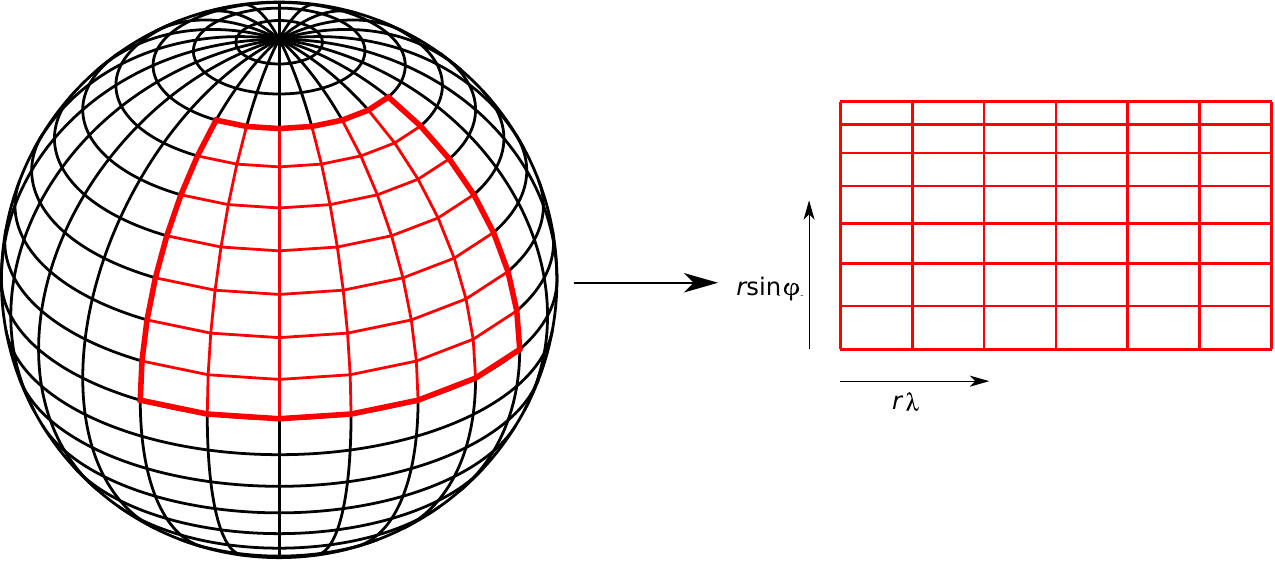}
	\caption{Schematic of transformation from $(\lambda,\varphi)$ to $(x,y)$ coordinate system for a long-lat region.}
	\label{fig:transformation}
\end{figure}

\paragraph{Forecasting results}
In order to propagate COD using $\u$, the following change of variable is made:
\begin{equation}
\label{eq:curvilinearVar}
\begin{split}
x &= r \lambda, \\
y &= r \sin \varphi,
\end{split}
\end{equation}
where $r$ is the radius of the earth, $\lambda$ denotes longitude and $\varphi$ denotes latitude. Figure~\ref{fig:latLongGlobe} shows the standard long-lat coordinate system, and Figure~\ref{fig:transformation} illustrates the coordinate transformation given by~\eqref{eq:curvilinearVar}.
The new variables, $x$ and $y$, are thus ``space'' variables corresponding to the longitude, $\lambda$, and latitude, $\varphi$, respectively. Note that at the poles, the transformation,~\eqref{eq:curvilinearVar}, is invalid since the longitude coordinates coalesce. Thus, it is not suitable for regions close to either pole.

The velocity field obtained through the optical flow procedure is in units pixel per time unit. In order to use this field for COD propagation, it must first undergo two conversions/transformations. First, we must convert the field to km/h. The time conversion is trivial as we know that the images are 15 minutes apart. For the space conversion, we must take into account how much space each pixel represents depending on its position. Assuming the earth to be spherical, the dimensions of a pixel do not depend on its longitudinal coordinate. However, its dimensions do depend on the latitudinal coordinate in the following way. The length of a degree of latitude does not vary much with the latitudinal coordinate, so we assume that each degree of latitude represents the same distance regardless of its coordinates. So for each pixel, we have
\begin{equation}
\label{eq:latLength}
\Delta x_{lat} = r\Delta \varphi
\end{equation}
where $\Delta x_{lat}$ is the distance a pixel represents along a line of longitude (i.e., distance along \emph{latitude} axis) and $\Delta \varphi$ is the latitude pixel spacing in radian measure which is uniform over the image.

The length of a degree of longitude, however, varies significantly with the latitude coordinate. We have
\begin{equation}
\label{eq:longLength}
\Delta x_{long} = \Delta \lambda \cdot r \cos \varphi
\end{equation}
where $\Delta x_{long}$ is the distance a pixel represents along a line of latitude and $\Delta \lambda$ is the longitude pixel spacing in radian measure which is also uniform over the image.

Having used the above the obtain the field, $\u$ in km/h, we must transform it from the $(\lambda,\varphi)$ coordinate system to the $(x,y)$ coordinate system. To do so, the following is applied:
\begin{equation}
\label{eq:velTrans}
u = U/ \cos \varphi, \, v = V \cos\varphi.
\end{equation}
The new velocity field, $(u,v)$ is valid on the $x$-$y$ grid. Note that both velocity components are transformed based on the latitudes of the velocity data on the original long-lat grid, and not on the longitudes. Looking at Figure~\ref{fig:transformation}, it is clear that moving away from the equator ($|\varphi|$ increasing), the transformation, \eqref{eq:curvilinearVar}, compresses the grid longitudinally, and expands it latitudinally. The velocity transformation, \eqref{eq:velTrans}, accounts for this by respectively increasing and decreasing the magnitudes of the latitudinal and longitudinal components of the velocity field with increasing $|\varphi|$ in transforming from the $(\lambda,\varphi)$ coordinate system to the $(x,y)$ system.

To advect COD images we employed the algorithm of Section~\ref{sec:dG}. We set $\Omega:=[-15567,-13788]\times[4009,4951]$ and discretize $\Omega$ with $200\times 200$ elements of order $N=3$. The timestep was selected depending on the advection coefficients speed to ensure that the CFL condition holds. The advection starts at 22:15 and uses the predicted velocity field $\u$ from NSE to predict dynamics of COD images up to 23:30. We compared this prediction against two most popular predictions, namely the so called persistency forecast (the latest available COD image is used for the forecast) and optical flow forecast (map $\F_4$ displayed at Figure~\ref{fig:F4NSE}). The prediction results are summarized in the Table~1 which shows that the NSE prediction outperforms the aforementioned popular predictors.

\begin{table}
 \begin{center}
 \begin{tabular}{| l | l | l | l|}
 \hline
{\bf Date: 2013-09-01}     &      NSE           &      Persistence &      Optical flow         \\ \hline
22:30 	 &	\textbf{34} 	 &	\textbf{57} &      \textbf{38}   	\\\hline
22:45 	 &	\textbf{43} 	 &	\textbf{71} &      \textbf{49}	 	\\\hline
23:00 	 &	\textbf{42} 	 &	\textbf{77} &      \textbf{48}	 	\\\hline
23:15 	 &	\textbf{48} 	 &	\textbf{100} &      \textbf{56}	 	\\\hline
\end{tabular}
\caption{Relative mean absolute percentage errors: ${100\%\sum_{s\in I}|C_{\mathrm{true}}(\bm{x}_s,t)-C_{\mathrm{predicted}}(\bm{x}_s,t)|}\,/\,{\sum_{s\in I}|C_{\mathrm{true}}(\bm{x}_s,t)|}$.}
\end{center}
\end{table}
Here $I$ is a set of indices of the grid points, which exclude the outliers. As noted above, in Section~\ref{sec:optflow}, some COD image values are abnormally high due to acquisition problems. These points (about 10\%) are therefore excluded from the assessment.

\section{Conclusion}
\label{sec:conclusion}

This work presents a new algorithm for dynamic cloud optical density forecast from satellite images. The algorithm relies upon ideas from data assimilation, numerical analysis and classical image processing. It shows good performance on real satellite images and outperforms the popular forecasting algorithms in the considered test cases. A challenging topic for the future research is to develop a more accurate model for the velocity field, e.g. take into account the fact that the 2D projection of a 3D velocity field could be compressible, and to improve upon the linear advection model, e.g. to introduce uncertain sources/sinks allowing the COD images changing intensity due to drops/increases in temperature.

\paragraph{Acknowledgement} 
This research is partially supported\footnote{\tiny This report was prepared as an account of work sponsored by an agency of the United States government. Neither the United States government nor any agency thereof, nor any of their employees, makes any warranty, expressed or implied, or assumes any legal liability or responsibility for the accuracy, completeness, or usefulness of any information, apparatus, product, or process disclosed, or represented that its use would not infringe privately owned rights. Reference herein to any specific commercial product, process, or service by trade name, trademark, manufacturer, or otherwise does not necessarily constitute or imply its endorsement, recommendation, or favoring by the United States government or any agency thereof. The views and opinions of authors expressed herein do not necessarily state or reflect those of the United States government or any agency thereof.} by US Department of Energy Contract NDE-EE0006017.

\appendix

\section{Gradient of $J$}
\label{sec:gradJ}

Let $\hat{\bm\omega}:=\{\hat\omega_{pq}\}_{|p|\le\frac{N_x}2,|q|\le\frac{N_y}2}$ be the vector of the projection coefficients of $P_{N_x,N_y}\hat\omega$, and define $\hat{\bm q}$, $\bm{v}_k$, $\bm{u}_k$ analogously. Note that $\hat{\bm\omega} = \Psi_{\frac1{\sqrt{2}},1} \Psi_{\frac1{\sqrt{2}},1}^\star\hat{\bm\omega}$, provided $\Psi_{a,b}$ is an $(N_x+1)(N_y+1)$ projection matrix defined by: $\Psi_{a,b}=
a\left( \begin{smallmatrix}
I&0&bJ\\0&0&0\\ J&0&-bI
\end{smallmatrix}
\right)$, i.e. $\hat{\omega}_{-\frac{N_x}2+k,-\frac{N_y}2+s}=\overline{\hat \omega}_{\frac{N_x}2-k,\frac{N_y}2-s}$ for any $0\le k,s\le \frac{N_x}2,\frac{N_y}2$ and $\hat{\omega}_{0,0}=0$ ($\hat\omega$ has zero mean!).
We derive the expression for the gradient of $J$ by using a so-called adjoint method: one computes Gateaux derivative of $J$ w.r.t. $\hat{\bm q}$, obtaines a linear ODE for the gradient, and, finally, uses its adjoint equation to find the gradient explicitly ($m:=(N_x+1)(N_y+1)$):
\[
\begin{split}
&J_{\hat u}(\hat{\bm q},\hat u,\hat v)=2\sum_{k=1}^{N_c} (d_k,\bm{\omega}^{\hat u}(t_k))_{C^m} - (u_k - \hat u,1)_{\Lt} \\
&\dfrac{d\bm{\omega}^{\hat u}}{dt} = -B(\left[\begin{smallmatrix}1\\0
\end{smallmatrix}\right])\hat{\bm\omega}-B_1(\hat{\bm\omega})\bm{\omega}^{\hat u} - B(\bm u) \bm{\omega}^{\hat u}- B(\hat{\bm\omega})\bm{\omega}^{\hat u} - A\bm{\omega}^{\hat u}\,,
\end{split}
\]
\[
\begin{split}
&J_{\hat v}(\hat{\bm q},\hat u,\hat v)=2\sum_{k=1}^{N_c} (d_k,\bm{\omega}^{\hat v}(t_k))_{C^m} -(v_k - \hat v,1)_{\Lt}\,, \\
&\dfrac{d\bm{\omega}^{\hat v}}{dt} = -B(\left[\begin{smallmatrix}0\\1
\end{smallmatrix}\right])\hat{\bm\omega}-B_1(\hat{\bm\omega})\bm{\omega}^{\hat v} - B(\bm u) \bm{\omega}^{\hat v}- B(\hat{\bm\omega})\bm{\omega}^{\hat v} - A\bm{\omega}^{\hat v}\,,\\
& \bm{\omega}^{\hat u} (0)=\bm{\omega}^{\hat v} (0)=0\,,
d_k:=A_x^\star(\bm{v}_k+A_x\hat{\bm\omega}(t_k)) - A_y^\star(\bm{u}_k-A_y\hat{\bm\omega}(t_k))\,,
\end{split}
\]
\[
\begin{split}
&\nabla_{\hat{\bm q}} J(\hat{\bm q},\hat u,\hat v)=2\sum_{k=1}^{N_c}\Psi_{1,\imath}^\star\bm{\lambda}_k(0)\,,\\
&\dfrac{d\hat{\bm \omega}}{dt} = -B(\hat{\bm \omega})\hat{\bm\omega} - B(\tilde u) \hat{\bm\omega} - A\hat{\bm \omega}\,, \hat{\bm\omega}(0)=\hat{\bm q}\,,
\\
  &\dfrac{d\bm{\lambda}_k}{dt} = B_1^\star(\hat{\bm\omega})\bm{\lambda}_k - B(\hat{\bm\omega})\bm{\lambda}_k - B(\tilde u) \bm{\lambda}_k + A\bm{\lambda}_k\,,\\
& \bm{\lambda}_k(t_k)=A_x^\star(\bm{v}_k + A_x\hat{\bm\omega}(t_k))-A_y^\star (\bm{u}_k - A_y\hat{\bm\omega}(t_k)\,,\\
&A_{x,y}:=\operatorname{diag}(\{\frac{2\pi i n}{\lambda_{nm}L_{x,y}}\}_{|n|\le \frac{N_x}{2},|m|\le\frac{N_y}{2}})\,,\\
&B_1(\hat{\bm\omega}) = [B(\Psi_{\frac1{\sqrt{2}},1} e_1) \hat{\bm\omega}\dots B(\Psi_{\frac1{\sqrt{2}},1} e_{m}) \hat{\bm\omega}]\Psi^\star_{\frac1{\sqrt{2}},1}\,.
\end{split}
\]
Note that all the differential equations above are linear, i.e. they can be represented in the form $\dot x = Q(t)x + s$, but the equation for $\hat{\bm \omega}$. This latter equation can be approximated by the following implicit time integrator: $\hat{\bm \omega}_{0} = \hat{\bm q}$ and
\begin{equation}
  \label{eq:NSE_midpoint}
\frac{\hat{\bm \omega}_{t+1}- \hat{\bm \omega}_{t}}{dt} = \left(-B(\hat{\bm \omega}_t) -B(\tilde u)-A \right) \frac{\hat{\bm \omega}_{t+1}+ \hat{\bm \omega}_{t}}{2}\,.
\end{equation}
$B(\hat{\bm \omega}_n)$ and $B(\tilde u) $ are skew-symmetric matrices, and $-A$ is a symmetric non-negative definite matrix, so that the above time integrator is unconditionally stable for $dt>0$, and it converges to the exact solution provided $dt\to0$~\cite{ZhukTTCDC15}. Once the discrete in time sequence $\{\hat{\bm \omega}_{t}\}$ is computed, one can approximate solutions of the aforementioned linear equations by using the standard simplectic midpoint method: \[
\frac{x_{t+1}- x_{t}}{dt} = \frac{Q_{t+1}+Q_t}{2}\frac{x_{t+1}+ x_{t}}{2}+\frac{f_{t+1}+f_t}2\,.
\]


\end{document}